\documentclass[authoryear,review]{elsarticle}
\usepackage{amsmath,amssymb}
\usepackage{makecell}
\usepackage{xcolor}
\usepackage{verbatim}

\usepackage{etoolbox}
\makeatletter
\patchcmd{\pprintMaketitle}
{\hrule\vskip12pt}
{\clearpage\hrule\vskip12\p@}
{}{}
\patchcmd{\pprintMaketitle}
{\hrule\vskip12pt}
{\hrule\clearpage}
{}{}
\makeatother

\journal{Expert Systems with Applications}

\biboptions{authoryear}

\begin{document}

\begin{frontmatter}

   \title{A multi-period and bi-objective approach for locating ambulances: a case study in Belo Horizonte, Brazil}
   \author[label0]{Charles Paulino de Oliveira} 
\ead{charles.paulino@outlook.com}
    \author[label1]{Elisangela Martins de S\'{a}*}
    \ead{elisangelamartins@cefetmg.br}
\author[label2]{Fl\'{a}vio Vin\'{i}cius Cruzeiro Martins} 
\ead{flaviocruzeiro@cefetmg.br}
\address[label0]{Postgraduate Program in Administration, Centro Federal de Educação Tecnológica de Minas Gerais, Av. Amazonas 7675, Belo Horizonte, MG, Brazil, 30510-000. \fnref{cor1}}
\address[label1]{Applied Social Sciences Department, Centro Federal de Educação Tecnológica de Minas Gerais, Av. Amazonas 7675, Belo Horizonte, MG, Brazil, 30510-000. \fnref{cor1}}
\address[label2]{Computer Department, Centro Federal de Educação Tecnológica de Minas Gerais, Av. Amazonas 7675, Belo Horizonte, MG, Brazil, 30510-000.}






\begin{abstract}
This work aims to apply the Facility Location Problem in the Emergency Medical Service (EMS) of Belo Horizonte, Brazil. The objective is to improve two previous optimization models from literature to handle base locations and ambulances allocation/relocation problems. The proposed multi-period models introduce the concept of relocation to the local EMS, which allows ambulances to move among bases in different periods to raise the system coverage. This paper also proposes a bi-objective approach aiming to minimize the number of bases and maximize the coverage of demands, which is solved using the $\epsilon$-constraint method.
Results show that coverage levels increase by up to 31\% for the deterministic approach and up to 24\% for the probabilistic approach. Ambulances' relocation optimization might improve coverage levels by up to 21\% for the deterministic approach and up to 15\% for the probabilistic approach when change scenarios from static to multi-period. Also, since the multi-period model solutions result in installing a larger number of bases, the bi-objective approach is a powerful tool for the decision-maker. Bi-objective results suggest modest increments in the objective function when the stations' number exceeds 28 for the deterministic approach. The probabilistic approach increments in the objective function start to narrow above 30 installed stations.
\end{abstract}

\begin{keyword}
Facility Location; Ambulance Location; Mathematical Programming Model; Real-World Application.
\end{keyword}

\end{frontmatter}


\section{Introduction}

Emergency Medical Service (EMS) is an essential component of healthcare systems. Its primary roles are to provide first aid to emergency medical demands out-of-hospital and transport victims to hospitals.
Some factors are essential to ensure successful treatment for those who need this kind of service: early detection, early reporting, early response, good on-scene care, care in transportation and transfer to definitive care. Among these factors, the early response is directly connected to the ambulances system's performance, which is sent to rescue the victims \citep{aringhieri2016, aboueljinane:13}.

  The response time indicator is commonly used to assess ambulances services' performance concerning efficiency in prehospital care. Response time can be accounted as the time interval between the moment the EMS receives an emergency call to the moment the patient is reached by ambulance. This kind of indicator is one of EMS systems' most critical indicators, especially in life-threatening emergencies. Studies about the impact of response times show that, on average, a minute increase in response times increases mortality by between 8 and 17\% - women and those over age 65 appear to be most affected \citep{wilde2013}.

  When planning an EMS, a predefined response time can be used to define the concept of coverage, where patients are said to be covered whether they can be reached by ambulance within the response time. In other words, coverage is defined as the proportion of demand that can be served within a given response time \citep{sudtachat2016}. Therefore, it is fundamental that ambulances be installed (allocated) in places that enable them to adequately cover their demands, making it possible to access demand points within sufficient response times.  In this sense, EMS managers worldwide face the difficult task of locating ambulance bases and allocating ambulances to the bases, which are limited resources to properly handle emergency medical demands.

Numerous studies have been carried out to achieve this goal, and a variety of tools like simulation, mathematical programming, and queuing theory models have been applied \citep{aboueljinane:13}. For decision-making, these applications reach various planning problems at the strategic, tactical, and operational level \citep{kergosien:15}. Even though the uncertain nature of demands makes the optimization of EMS systems a complex task, applying mathematical programming tools for the  EMS location problem has attracted particular attention in the last few decades.

Most mathematical programming researches focus on strategic and tactical issues, which are problems of a static nature that handle decisions of determining where ambulance bases should be located and the numbers of ambulances per base \citep{van2016, belanger2019}. Some other works address short-term issues decisions at the operational level, such as relocation strategies and vehicle dispatching \citep{van2016, belanger2019}. This paper focuses on EMS strategic, tactical, and operational levels of decision-making. The proposed approach concerns two decision problems simultaneously: installing the ambulance bases and how and when ambulances should be allocated to bases.

This work improves two previous literature approaches that address locating bases and allocating ambulances over bases problem. The improvement consists of adding the relocation strategy to these models. The relocation strategy implemented allows to consider the dynamism of demand in multiples periods, e.g., hourly in the day, and let ambulances change their locations at the end of a period to achieve higher coverage in the next period.

The two problems consider two different assumptions. The first approach is a deterministic coverage model that tries to maximize the coverage, assuming that a demand node is always covered when any bases/ambulance is placed on its coverage zone. Meanwhile, the second formulation considers the ambulance availability and seeks to maximize the expected coverage considering the ambulance's busyness likelihood.  

Furthermore, since both approaches do not restrict the number of bases installed, a multi-objective approach that aims to maximize the coverage of demand and minimize the number of bases is proposed for both models. The $\epsilon$-constraint method is applied to find the Pareto-optimal solutions.

Hence, the main contribution of the present study is twofold. First, the study presents the improvement of two static optimization models from the literature considering multiple periods and the optimal ambulances allocations in each period. The models take into account time-dependent demand patterns and also consider multiple types of ambulances with different coverage requirements. Second, the work proposes a bi-objective problem that aims to maximize demand coverage and minimize the number of ambulance bases. Further, the proposed models are analyzed using data from the EMS from Belo Horizonte city, Brazil.

The remainder of this paper is organized as follows. Section \ref{literature} shows a literature review of facility location problems for EMS. Section \ref{formulation} presents the notation and the two proposed formulations.
Section \ref{biobjective} presents the multi-objective models.
Section \ref{res} presents the computational experiment results and its analysis, confirming the relocation strategy's effectiveness.
The paper finishes with a section of final remarks and possible future research lines. 


\section{Facility Location Problem for Emergency Medical Services}
\label{literature}

Applications of the facility location problem include a diverse set of examples as EMS, fire stations, schools, hospitals, reserves for endangered species, airline hubs, waste disposal sites, and warehouses \citep{daskin:08}.
The facility location problems applied to EMS can be categorized into three class: 
(\textit{i}) covering models, which focus on demand points coverage within a predefined distance or time limit (response time); 
(\textit{ii}) $p$-median models, which aim to minimize the total or average  distance, cost, or service time for all demand points; and 
(\textit{iii}) $p$-center models, which minimize the maximum distance, costs, or travel time between opened facilities and demand points allocated to each one \citep{jia:07, li:11}.

Covering models have always been very attractive in EMS facilities problems due to its easy applicability in general real-world problems. The literature considers two types of covering models: set covering and maximal covering. As one of the early contributors in this research area, \citet{toregas1971} introduce the Set Covering Location Problem, which focuses on minimizing the number of ambulances needed to cover all demand points within a given distance or time threshold. \citet{church:74} propose a model, known as Maximal Coverage Location Problem (MCLP), which aims to maximize the demand points covered given a fixed-sized fleet. 

A significant number of variants from these two seminal models have been proposed to address EMS location problems. The MCLP approach has been widely used in practice, and it may be the most influential of all EMS location models \citep{erkut2009}. Therefore, this work concentrate on its branches. Reviewer papers \citep{brotcorne:03, li:11, belanger2019} show the modeling evolution and give a detailed literature overview. \citet{brotcorne:03} suggest categorizing these models according to how data stochasticity is incorporated into the model. They could be deterministic or probabilistic. They also propose a classification of these models as static or dynamic, where static models address base location or vehicle allocation decisions only once while dynamic models allow different base locations or vehicle allocations in different periods.

Most models for EMS location problems belong to the class of static problems, where decisions are taken to establish the optimal location to stations, assign ambulances to stations, and determine the fleet size \citep{aringhieri2017, belanger2019}. According to the location plan, static models assume that each ambulance is sent to its predefined station whenever it becomes idle. Nonetheless, since the demands vary spatially and temporally by day of the week, and by the time of the day,  to maximize the coverage of emergency demands,  in a given period, idle ambulances assigned to low demand areas may need to move to busier areas \citep{li:11}.

Hence, dynamic models have been developed to relocate ambulances throughout the day to consider the system's evolution over time \citep{brotcorne:03}. Dynamic models can be divided into two categories: multi-period and real-time \citep{van2018}. The multi-period model consists of establishing different location plans for each period, in which ambulances can be scheduled to move between periods to implement the next location plan. The real-time model consists of relocating available ambulances whenever one or more vehicles are dispatched to respond to an emergency demand. In this case, a relocation decision depends on the system's state, in real-time, and relocation may be needed, at any time, to sustain a proper service level.

Most of the early models do not consider the stochastic nature of EMS operations, ignoring the availability of the ambulance \citep{brotcorne:03}. These models are classified as deterministic and assume that ambulances are always available to be dispatched upon the arrival of emergency demands. However, this is not always the case in practice since vehicles' assignment to other emergency demands may let initially covered areas uncovered.

Based on a deterministic and static approach, \citet{schilling:79} propose two important models, the Tandem Equipment Allocation Model (TEAM) and the Facility-Location Equipment-Emplacement Technique (FLEET). Both models consider a fixed number of vehicles of different types, which have to be allocated to fire stations. These models are initially developed to allocate different vehicle types to fire stations, but they are also relevant and well known in the EMS context. The FLEET formulation also assumes a predefined number of stations to be installed; since the opening and conservation of stations might involve costs, it is appropriate to fix the stations' number. The first model of this paper is based on the FLEET model because of its ability to consider different equipment and limitations on their quantity, standard in EMS systems with different types of ambulances.

To represent the FLEET formulation, $I$ denotes the set of demand points, and $J$ denotes the set of potential ambulance location sites. The sets $N_i^p \subset J $ and $N_i^s \subset J $ correspond to the subsets of location sites that can cover a demand point $i$ $\in$ $I$ with vehicles type $p$ and $s,$ respectively. In this case, a demand point $i \in I$ is said to be covered by any location site $j \in J$ if and only if the distance between them, given by $r_{ji},$ is within the coverage limit represented by response time parameters $S^p$ and $S^s$ for vehicle type $p$ and $s,$ respectively, i.e, if $r_{ji}$ $\leq$ $S^p$ and $r_{ji}$ $\leq$ $S^s$.

The number of vehicles type $s$ and $p$ to be allocated to stations and the number of stations to be installed are represented by $P^s$, $P^p$, and $P^z$, respectively. The demand associated with each demand point $i$ $\in$ $I$ is denoted by $d_i$. The binary variables $x_{j}^p$  and $x_{j}^s$ is equal to 1 if vehicle type $p$ and $s$, respectively, is located at node $j \in J$ and 0, otherwise. The binary variable $y_i$ is equal to 1 if the demand point $i \in I$ is covered by both vehicle types ($s$ and $p$) and 0, otherwise. Finally, the binary variable $z_j$ is equal to 1 if a facility is located at the location site $j \in J$ and 0, otherwise. The FLEET model, propose by  \cite{schilling:79}, can be formulated as follows:

\begin{align}
\text{Maximize $F =$} \ & \sum_{i \in I} d_{i} y_{i} \label{foF} \\
\label{rF1} \mbox{subject to:} \ & \sum_{j \in N_{i}^p} x_{j}^p \geq y_{i} &&  \forall   i \in I \\
\label{rF2} \ & \sum_{j \in N_{i}^s} x_{j}^s \geq y_{i} && \forall   i \in I  \\
\label{rF3} \ & \sum_{j \in J} x_{j}^p = P^p  \\
\label{rF4} \ & \sum_{j \in J} x_{j}^s = P^s  \\
\label{rF5} \ & \sum_{j \in J} z_{j} = P^z  \\
\label{rF6} \ & x_{j}^p \leq z_{j}&& \forall j \in J\\
\label{rF7} \ & x_{j}^s \leq z_{j}&&  \forall j \in J \\
\label{rF8} \ & x_{j}^s, x_{j}^p \in \{0, 1\} &&  \forall j \in J \\
\label{rF9} \ & y_{i} \in \{0, 1\} &&  \forall i \in I \\
\label{rF10} \ & z_{j} \in \{0, 1\} &&  \forall j \in J.
\end{align}

The objective function \eqref{foF} maximizes the sum of covered demands. The first and second set of constraints, \eqref{rF1} and \eqref{rF2}, ensure that a demand point can only be considered covered if it is covered by both vehicle types $p$ and $s$. Constraints \eqref{rF3}, \eqref{rF4} and \eqref{rF5} limit the number of vehicles ($p$ and $s$) to be allocated and stations  to be installed. The set of constraints \eqref{rF6} and \eqref{rF7} prohibit the emplacement of equipment at location sites where stations have not been located. Constraints \eqref{rF8}, \eqref{rF9} and \eqref{rF10} set the variables as binary variables. 
 
TEAM and FLEET are static and deterministic location models. Because of that, they have a shortcoming of adequately account for the coverage when ambulances become busy. Two strategies emerge to address this problem: one providing multiple coverages and the other, taking into account the busyness probabilities of facilities/ambulances \citep{li:11, belanger2019}. Multiple coverage models are also a deterministic approach. However, they seek to increase the likelihood of having a demand zone covered by one available ambulance by increasing the number of ambulances covering this zone \citep{belanger2019}. The second strategy refers to the class of probabilistic models based on the busy fraction of the vehicles, which are the fraction of time that an ambulance is not available to answer demand and can be estimated in several ways \citep{brotcorne:03}.

\citet{gendreau:97} introduce the Double Standard Model (DSM). It is one of the most important models addressing multiple coverages. The DSM aims to maximize the demand covered twice within a response time goal $S_1$, establishing that at least a proportion of the demand must be covered within this response time $S_1$. A less strict response time $S_2$ is also set and all the demand must be covered within this secondary level of coverage ($S_2$ $>$ $S_1$). As evolution of DSM, \citet{schmid2010} present a multi-period version of the DSM and \citet{gendreau:01} present a real-time approach.

Probabilistic models are an alternative to approach real-life situations of ambulance unavailability. This strategy includes the Maximum Expected Covering Problem (MEXCLP) propose by \cite{daskin:83}, and the Maximum Availability Location Problem (MALP) from \cite{revelle1989}. Whereas the MEXCLP model seeks to maximize the expected coverage given that servers (ambulances) might be busy, the MALP model aims to maximize the sum of demands that can be covered with some minimum level of reliability ($\alpha$), also considering that ambulances might be busy \citep{sorensen2010}. Either MEXCLP or MALP might request to allocate more than one ambulance to cover demand points, depending on its demand level and the number of vehicles available.

Due to the model's probabilistic nature, stochastic elements of the EMS systems result in non-linear mathematical models. Consequently, this becomes the main difficulty. Hence, some assumptions might be necessary to incorporate probabilistic information in linear programming models and model its behavior with linear expressions \citep{sorensen2010}. The most common assumptions are: 
$(i)$ the ambulances busyness is uniform, i.e., all vehicles are equally busy, despite their location;  
$(ii)$ the ambulances are independent, which means that the probability that a vehicle is available is independent of the availability of other vehicles; and 
$(iii)$ the service areas are locally-constrained, i.e., ambulances within a local area only handle demands within that same local area \citep{sorensen2010}.

Simplifying assumptions regarding the uniform server busyness and the independence of the servers are made in MEXCLP. The MALP model also brings the assumption of server independence and includes the locally constrained service area assumption, but it relaxes the assumption of uniform ambulance availability. In MALP, \citet{revelle1989} introduce a method to compute an estimate of area-specific busy fractions associated with each demand point. It is possible to consider the ambulance availability in the coverage area of each demand node.

From the MALP and MEXCLP framework, some other models have emerged to improve the probabilistic approach. Some relevant examples are the adjusted MEXCLP model (AMEXCLP) from \cite{batta1989maximal}, the Queuing Maximal Availability Location Problem (Q-MALP) proposed by \cite{marianov1996}, the Local Reliability-MEXCLP (LR-MEXCLP) presented by \cite{sorensen2010}, the MEXCLP model with time variation (TIMEXCLP) formulated by \cite{repede:94}, the dynamic MEXCLP (DMEXCLP) introduced by \cite{jagtenberg2015} and The Facility Location and Equipment Emplacement Technique model with Expected Covering (FLEET-EXC) from \cite{rodriguez2020} which also improves the FLEET model.

The LR-MEXCLP model integrates the area-specific busy fraction of MALP with the maximum expected coverage objective of MEXCLP. Since the most common EMS goal is to maximize aggregate system response across all demands, \cite{sorensen2010} argue that the MALP's $\alpha$-reliability objective function may not be appropriate in some EMS contexts and may lead to inferior results. However, MALP's local reliability estimates constitute a useful modeling tool and outperform the uniform server busyness assumption of MEXCLP. For this reason, the LR-MEXCLP is used to incorporate stochasticity into the FLEET model, giving rise to the second problem addressed in this article.

As proposed by \cite{sorensen2010}, to model the LR-MEXCLP,  let the neighborhood $\bar{N}_i \subseteq I$ of the demand point $i \in I$ denote the set of demand points $j \in I$ where the travel time from $i$ to $j$ is less than or equal to the target response time, i.e., $r_{ji} \leq S$.  Also, let $\bar{d}_j$ represent the amount of demand at demand point $j \in I$ measured in service time. Hence, the area-specific busy fraction of demand point $i$ when $k$ servers are allocated within $\bar{N}_i$ is represented by $b_{ik},$ which can be computed by means of expression \eqref{rAA}. Further, let $q_{ik}$ denote the likely service reliability  at demand point $i \in I$ given the location of $k$ servers within its neighborhood $\bar{N_i},$ which can be computed by means of expression \eqref{rAB}.

\begin{align}
\text{$b_{i, k}$ = } \frac{\sum_{j \in \bar{N}_i} \bar{d}_j}{k} \label{rAA}
\end{align}

\begin{align}
\text{$q_{i, k}$ = } 1 - (b_{i, k})^k.   \label{rAB}
\end{align}

In LR-MEXCLP, as in the FLEET model, $d_i$ expresses the amount of demand at each demand point $i$ $\in$ $I$. Through the formulation, $N_i \subseteq J$ indicate the subset of location sites that can cover a demand point $i$ $\in$ $I$ within the coverage standard $S$ ($r_{ji}$ $\leq$ $S$). Let $P$ denote the maximum number of ambulances to be located, and let $K=\{1,\ldots,P\}$ denote the possible numbers of ambulances that can cover a demand point. The binary variable $x_{j}$ is equal to 1 if an ambulance is located at node $j$ $\in$ $J$ and 0, otherwise. Yet, the binary variable $y_{ik}$ is equal to 1 if demand point $i$ $\in$ $I$ is covered by $k$ servers and 0, otherwise. The LR-MEXCLP model, proposed by \cite{sorensen2010}, is formulated as:

\begin{align}
\text{Maximize $F =$} & \sum_{i \in I} \sum_{k \in K} d_{i} q_{ik} y_{ik}  \label{fL}\\
\label{rL1} \mbox{subject to:}  \ & \sum_{j \in N_{i}} x_{j} - \sum_{k \in K} k y_{ik} \geq  0 && \forall i \in I  \\
\label{rL2} \ & \sum_{k = 1}^{P} y_{ik} \leq 1 && \forall i \in I \\
\label{rL3} \ & \sum_{j \in J} x_{j} \leq P &&  \\
\label{rL4} \ & x_{j} \in \{0, 1\} && \forall j \in J \\
\label{rL5} \ & y_{ik} \in \{0, 1\} &&\forall i \in I, \forall k=1,\ldots,P.
\end{align}

The expression \eqref{fL} aims to maximize the sum of demands multiplied by the likely reliability of service. As pointed out by \cite {sorensen2010}, the objective function estimates the amount of demand that will be met for a  target response time. The first set of inequalities \eqref{rL1} is used to compute the number of ambulances by which a node is covered. In this case, if the sum of the $x_{j}$ values is less than $k$, then it will not be feasible to set the corresponding $y_{ik}$ value to 1. Given that $q_{ik}$ increases with the increment of $k$ values and the maximization objective, $y_{ik}$ will always attempt to set 1 for the largest $k$ value possible for each node $i$. The set of inequalities \eqref{rL2} establish that at most one $y_{ik}$ value, for any given node $i$, can be set to 1. The expression \eqref{rL3} restricts the employment of servers to its limit. Finally, constraints \eqref{rL4} and \eqref{rL5} set the variables as binary variables. 

Although LR-MEXCLP relaxes the uniform server busyness assumption, it still relying on assumptions of server independence and locally-constrained service areas. Nevertheless, the LR-MEXCLP model offers improvements over both MEXCLP and MALP. Further, as stated by \cite{aringhieri2017}, it is difficult to determine which simplified assumption imposes more limitations on finding more realistic solutions. As correctly mentioned by \cite{sorensen2010}, despite efforts to improve probabilistic models, relaxing one or more assumptions often leads to another imposition.


\section{Problems description and mathematical models}
\label{formulation}

This paper propose improvements for two EMS location/allocation models: the FLEET model \citep{schilling:79} and the LR-MEXCLP model \citep{sorensen2010}. The objective is to give the population the best EMS possible, maximizing the coverage, assuming a limited number of ambulances to allocate on bases. Thus, it is intended to add the relocation strategy to both models, allowing them to meet the strategic, tactical, and operational levels of decision making in EMS systems. Besides, the concept of independent coverage in these models is also included.

The proposal for adding the relocation strategy to the formulations makes it possible to change ambulances allocation between periods, contrary to static models, which assume that ambulances are assigned to a home base and must return to that same base whenever a service is finished. By allowing relocations, it intends to consider demand variation through different periods, besides the variation of demand in distinct areas.

It is assuming that EMS systems operate with a set of service types, denoted by $U$, where each service type has its set of vehicles, response time standards, and demand patterns. Basing on the concept of independent coverage,  presented in   \cite{guimaraes2020}, it considers that the demand coverage of different service types is independent. So, a demand point may be covered by one type of ambulance without the obligation of being covered by the other type of ambulance. Moreover, this work assumes that an ambulance of a type $u \in U$ can only handle the demands of type $u$ service. An illustration is the EMS systems addressed in \cite{guimaraes2020} that operate with two types of servers, advanced life support (ALS) units, and basic life support (BLS) units. So, ALS and BLS units provide different services, and the services provided by each one require different response times.

Unlike \cite{repede:94}, \cite{schmid2010} and \cite{van2015}, the proposed approach in this paper do not consider the travel time variation between time periods. Instead, it assumes that travel times are deterministic. Moreover, the proposed models consider that once a base is opened, it must remain open until the end of the planning horizon and that there is no limit to the number of bases. Considering the situations that may not be convenient to open many bases, the minimization of the number of bases will be addressed in the multi-objective models presented in Section \ref{biobjective}. It also assumes a fixed number of ambulances of each type, and each base can hold a maximum number of vehicles.

The necessary basic notation not yet provided is described. However, the notation specific to each model is described along with its corresponding formulation. Let  $T$ denote the set of time periods. The number of each ambulances type $u$ $\in$ $U$ available to be allocated to a base is  $P_{u}.$ For each ambulance type $u$ $\in$ $U$, let $S_{u}$ represent its standard response time. The maximum number of vehicles that each base $j \in J$ can hold is $C_{j}.$ The  number of demands generated at each demand point $i$ $\in$ $I$ for each ambulance type $u \in U$ at each time period $t \in T$ is denoted by $d_{iu}^t$.

\subsection{Facility Location and Equipment Emplacement Technique with Independent Coverages and Time-Varying Demands (FLEET-ICt)}

The multi-period model FLEET-ICt (Facility Location and Equipment Emplacement Technique with Independent Coverage and Time-Varying Demands) is proposed based on the FLEET model and independent coverage concepts. The FLEET-ICt formulation has a deterministic approach, as in FLEET. It does not account for the ambulance's unavailability. It is assumed that ambulances are always available to handle their demands.

In FLEET-ICt, the binary variable $x_{ju}^t$ is equal to 1 if an ambulance type $u$ $\in$ $U$ is allocated to site $j$ $\in$ $J$ at time period $t$ $\in$ $T$ and 0 otherwise. The binary variable $y_{iu}^t$ is equal to 1 if a demand point $i$ $\in$ $I$ is covered by ambulance $u$ $\in$ $U$ at time period $t$ $\in$ $T$. Finally, the binary variable $j$ is equal to 1 if a station is installed at site $j$ $\in$ $J$ and 0, otherwise. The FLEET-ICt model can be formulated as follows:

\begin{align}
\text{Maximize $F_1$ =} & \sum_{i \in I} \sum_{u \in U} \sum_{t \in T} d_{iu}^t y_{iu}^t  \label{fo1}\\
\mbox{subject to:} \label{rA11} \ & \sum_{j \in N_{i}^u} x_{ju}^t \geq y_{iu}^t && \forall i \in I,   u \in U,   t \in T \\
\label{rA12} \ & \sum_{j \in J} x_{ju}^t \leq  P_{u} && \forall   u \in U,    t \in T  \\
\label{rA14} \ & \sum_{u \in U} x_{ju}^t \leq  C_{j} z_{j} && \forall   j \in J,   t \in T  \\
\label{rA15} \ & x_{ju}^t \leq z_{j} && \forall j \in J,   u \in U,  t \in T \\
\label{rA16} \ & z_{j} \in \{0, 1\}&&\forall j \in J \\
\label{rA17} \ & x_{ju}^t \in \{0, 1\}&&\forall j \in J,   u \in U,   t \in T \\
\label{rA18} \ & y_{iu}^t \in \{0, 1\}&&\forall i \in I,  u \in U,  t \in T.
\end{align}

The expression \eqref{fo1} aims to maximize the amount of demand covered, for each type of ambulance, for each period. The set of inequalities \eqref{rA11} ensure that a demand for an ambulance type can only be considered covered, in a given period, if at least one ambulance of this same type is allocated within the maximum response time standard. The inequalities set \eqref{rA12} restrict the employment of ambulances, for each type of ambulance, for each period to its limit. The inequalities set \eqref{rA14} limits the number of ambulances for each station for each period. The constraints \eqref{rA15} ensure that ambulances will not be allocated to points where there is no active station. Constraints \eqref{rA16}, \eqref{rA17} and \eqref{rA18} represent binary requirements.

\subsection{Local Reliability Maximum Expected Covering Problem with Independent Coverages and Time-Varying Demands (LR-MEXCLP-ICt)}

The Local Reliability Maximum Expected Covering Problem with Independent Coverage and Time-Varying Demands (LR-MEXCLP-ICt) is a multi-period extension of the LR-MEXCLP model, introduced by \cite{sorensen2010}. It also includes the FLEET approach of considering the allocation of different ambulance types and the concept of independent coverage.  

Let $K=\{1,\ldots, k_{max}\}$ represents the possible numbers of ambulances that can cover a demand point, where $k_{max}$ is the maximum number of ambulances that can cover a demand point.  Now the likely service reliability for ambulance type $u \in U$ at demand point $i \in I$ at time period $t \in T$ given the allocation of $k \in K $ ambulances within its neighborhood is designated as $q_{iu k}^t$.  The area-specific busy fraction of ambulance type $u \in U$ of those $k \in K$ ambulances within the same area surrounding $i$ at each period $t$ is then represented by $b_{iu k}^t$. The calculation of $b_{iu, k}^t$ and $q_{iuk}^t$ is given by \eqref{rAC} and \eqref{rAD}, where $\bar{d}_{ju}^t$ represent the amount of demand measured in service time for ambulance type $u$ at demand point $j \in I$ in period $t \in T.$

\begin{align}
\text{$b_{iuk}^t$ = } \frac{\sum_{j \in \bar{N}_i} \bar{d}_{ju}^t}{k} \label{rAC}   
\end{align}

\begin{align}
\text{$q_{iuk}^t$ = } 1 - (b_{iu, k}^t)^k   \label{rAD}
\end{align}

In LR-MEXCLP-ICt, the binary variable $x_{ju}^t$ is equal to 1 if an ambulance type $u$ $\in$ $U$ is allocated to site $j$ $\in$ $J$ at time period $t$ $\in$ $T$ and 0, otherwise. The binary variable $y_{iu k}^t$ is equal to 1 if a demand point $i$ $\in$ $I$ is covered by  $k \in K$ ambulances of type  $u \in U$  at period $t \in T$ and 0, otherwise. Finally, the binary variable $z_j$ is equal to 1 if a station is installed at site $j$ $\in$ $J$ and 0, otherwise. The new version of the LR-MEXCLP, called LR-MEXCLP-ICt is formulated as follows:

\begin{align}
\text{Maximize $F_2$ =} & \sum_{i \in I} \sum_{u \in U} \sum_{t \in T} \sum_{k \in K} d_{iu}^t q_{iuk}^t y_{iuk}^t  \label{fo7}\\
\label{rA36} \mbox{subject to:}  \ & \sum_{j \in N_{i}^u} x_{ju}^t - \sum_{k \in K} k y_{iuk}^t \geq  0 && \forall i \in I, u \in U ,   t \in T,  \\
\label{rA37} \ & \sum_{k  \in K} y_{iuk}^t \leq 1 && \forall i \in I,  u \in U,   t \in T  \\
\label{rA38} \ & \sum_{j \in J} x_{ju}^t \leq P_{u} && \forall u \in U,   t \in T  \\
\label{rA39} \ & \sum_{u \in U}  x_{ju}^t \leq C_{j} z_{j} && \forall j \in J,   t \in T \\
\label{rA41} \ & x_{ju}^t \leq z_{j} && \forall j \in J,  u \in U, t \in T \\
\label{rA42} \ & z_{j} \in \{0, 1\} && \forall j \in J \\
\label{rA43} \ & x_{ju}^t \in \{0, 1\} && \forall j \in J, u \in U,  t \in T \\
\label{rA44} \ & y_{iu, k}^t \in \{0, 1\} && \forall i \in I,  u \in u, k \in K, t \in T.
\end{align}

The expression \eqref{fo7} aims to maximize the sum, for all nodes, of the demand level weighted by the reliability of coverage. The first set of inequalities \eqref{rA36} is used to determine the number of servers by which a node is covered. In these constraints, if the sum of the $x_{ju}^t$ values is less than $k$, then it will not be feasible to set the corresponding $y_{iuk}^t$ value to 1. Given that the value of $q_{iuk}^t$ increases with $k$, the $y_{iuk}^t$ variable with the largest possible value of $k$ will always be set to 1 for each node $i$. By the set of inequalities \eqref{rA37}, it is established that at most one $y_{iuk}^t$ value, for any given node $i$, can be set to 1. The set of inequalities \eqref{rA38} restricts ambulances' employment for each period to its limit. The number of ambulances for each station in each period is limited by the inequalities set \eqref{rA39}.  Constraints \eqref{rA41} ensure that ambulances will not be allocated to points where there is no active station. Finally, Constraints \eqref{rA42}, \eqref{rA43} and \eqref{rA44} represent binary requirements.

\subsection{Comparison of FLEET-ICt and LR-MEXCLP-ICt problems}

The difference between FLEET-ICt and LR-MEXCLP-ICt lies on the coverage concept associated with each objective function and also on the employment of the $y$ variable. Table \ref{tab:my_label} summarizes a comparison between both models.

In the FLEET-ICt model, the objective function aims to maximize coverage through a deterministic approach. In this model, the coverage concept is only associated with an ambulance allocation (or not) within a coverage zone. Hence, the role of $y$ variable is to indicate if a demand point is covered (or not). 

On the other hand, the objective function in LR-MEXCLP-ICt considers that, due to the occupancy rate of the ambulances, the number of ambulances in the coverage zone of a demand point must be considered when estimating the amount of demand that will be served in a given response time. In this case, the objective function considers the amount of demand expected to be attended at a target service level, considering the ambulances' busy fraction. Here, the $y$ variable has the role of indicating if a demand point is covered (or not) and which service level (ambulance number - $k$) is assigned to each of them.

\begin{table}[h]
    \centering
    \begin{tabular}{p{1.8cm}|p{4.5cm}|p{4.5cm}} \hline
 Formulation & FLEET-ICt & LR-MEXCLP-ICt\\ \hline
        Objective function &  Maximize coverage: amount of demand covered & Maximize coverage: amount of demand covered   multiplying by relibiality of the service \\ \hline
        Coverage & Deterministic: take into account the  ambulance location into a  coverage zone ($r_{ji}$ $\leq$ $S_u$) & Probabilistic: take into account the ambulance location into a  coverage zone ($r_{ji}$ $\leq$ $S_u$) and the ambulance busy fraction  within a coverage zone \\ \hline
        Variable & $y_{iu}^t$: equal to 1 if a demand point  $i$ is covered by ambulance type $u$  at time period $t$ & $y_{iu, k}^t$: equal to 1 if a demand point $i$ is covered by ambulance type $u$  with $k$ vehicles at time period $t$ \\ \hline
    \end{tabular}
    \caption{Comparison of FLEET-ICt and LR-MEXCLP-ICt problems.}
    \label{tab:my_label}
\end{table}

\section{Bi-objective problem} 
\label{biobjective}

Real-world problems are inherently composed of several criteria that need to be optimized. These different criteria are placed as objective functions that are, in most cases, conflicting. Solutions in these cases often have multiple dimensions and must be analyzed despite their conflicts. Some examples of conflicting criteria are costs and time, productivity and delays, need and availability, coverage, and the number of resources.

Addressing an EMS context, \cite{schmid2010} propose a multi-objective formulation that aims to maximize coverage and, at the same time, must minimize the number of relocations, avoiding excessive ambulance changes through stations. 
\cite{su2015cost} address a multi-objective model that focuses on maximizing coverage and minimizing the cost attributed to delays in emergency attendance (considering a standard response time) and the operational cost of the EMS. \cite{guimaraes:18m} present a formulation aiming to maximize the number of demands covered and minimize the standard deviation of the average demand response time. 

A multi-objective approach for an EMS is addressed by \cite{ferrari2018} which propose a multi-objective formulation considering five objectives: maximizing the covered population, maximizing the number of demand points covered, minimizing the number of stations to be installed, minimizing the distance between demand points and stations, and also minimize the penalization for the demand that can not be allocated to any bases due to capacity restrictions. \cite{el2019bi} propose a bi-objective approach that maximizes the expected coverage while minimizes the coverage costs incurred when demands are assigned to another service provider (private company).  \cite{guimaraes2020} address a bi-objective approach that maximizes the covered calls and minimizes the number of bases installed.

In these cases, solutions might bear trade-offs between the different objectives. 
Therefore, a solution set, called Pareto, is required to represent the different optimal solutions (non-dominated solutions) considering all objectives of a multi-objective problem. As \cite{lobato2017} mention, extreme solutions on the Pareto set cannot satisfy all objective functions, and the optimal solution for one of the objectives will not necessarily be the best solution for the others.

\subsection{Bi-objective models}

This section considers the situations where an EMS operation must avoid excessive resource allocation. Therefore, the benefits of an increased number of covered calls have to be balanced with the cost of increase the system's capacity. Among these resources, the number of ambulances and the number of stations are essential.

According to the EMS resources limit, both formulations consider restrictions in employment regarding the number of ambulances. However, there is no parameter determining the maximum number of stations to be installed. The models allow the assignment of different stations for each ambulance in each period. Thus, with an increased number of periods, a more significant number of stations installed can be expected. A multi-objective approach is presented to avoid such situations. A second objective function \eqref{fo91} is add to the FLEET-ICt formulation \eqref{fo1}-\eqref{rA18} and to the LR-MEXCLP-ICt formulation \eqref{fo7}-\eqref{rA44}. This new objective function aims to minimize the number of stations to be installed.

\begin{align}
\textrm{Minimize $F_3$ = } & \sum_{j \in J} z_{j} \label{fo91}
\end{align}

Through the $\epsilon$-constraint method \citep{yv1971}, the second objective function is placed as a constraint \eqref{rA131}, which turns the bi-objective problem into a single objective problem. Here, $\epsilon$ indicates the number of stations that can be installed on each iteration.

\begin{align}
\label{rA131} \ & F_3 \leq \epsilon
\end{align}

Therefore, the scalarized single objective version of the bi-objective model for FLEET-ICt is given as:

\begin{align}
 \text{Maximize}   & \ \sum_{i \in I} \sum_{u \in U} \sum_{t \in T} d_{iu}^t y_{iu}^t \label{fo1001} \\
		\mbox{subject to:}  \ & \eqref{rA11}-\eqref{rA18} \nonumber\\
 \label{rA132} \ &    \sum_{j \in J} z_{j} \leq \epsilon.
 \end{align}

The objective function \eqref{fo1001} and the constraints \eqref{rA11} - \eqref{rA18}  work in the same way described on FLEET-ICt formulation. Constraint \eqref{rA132} limit the number of stations to be installed to $\epsilon$, according to the second objective function \eqref{fo91}.

Straightaway, the scalarized single objective model for LR-MEXCLP-ICt is given as:

\begin{align}
 \text{Maximize}   & \sum_{i \in I} \sum_{u \in U} \sum_{t \in T} \sum_{k \in K} d_{iu}^t q_{iu, k}^t y_{iu, k}^t \label{fo1000} \\
		\mbox{subject to:}  \ & \eqref{rA11}-\eqref{rA18} \nonumber\\
 \label{rA133} \ &    \sum_{j \in J} z_{j} \leq \epsilon.
 \end{align}

The objective function \eqref{fo1000} and the constraints \eqref{rA36} - \eqref{rA44} are the identically as presented before on LR-MEXCLP-ICt formulation. Finally, constraint \eqref{rA133} limit the number of stations to be installed to $\epsilon$, according to the second objective function \eqref{fo91}.

\section{Computational experiment results}
\label{res}

\subsection{Belo Horizonte's EMS dataset}

The database used in this paper came from   \citet{guimaraes2020}. It is composed of real data from Belo Horizonte's EMS, with service types composed of ALS and BLS units. The results consider the distribution of demand by different periods of the day differently from \citet{guimaraes2020} which address a static model and, consequently, consider only the geographic distribution of demand. This dataset has demands ($d_{iu}^t$) from a period range of one year. \citet{guimaraes2020} present all detail about the instance collected.

The set of demand nodes $I$ is composed of 427 Belo Horizonte's districts. The 427 demand points mapped are represented by the geometric center of each one of the considered districts.  The set of potential ambulance location sites $J$ is defined based on a search for buildings in the city with an adequate structure to install a station. Among the potential stations, there are hospitals, shopping centers, schools, parking lots, public buildings, etc. Thus, 1527 potential ambulance location sites are mapped as a possible choice to be a base station.

The set of periods $T$ considered is 24, composed of 1 hour each. The set of ambulance types $U=\{ALS, BLS\}$. The number of vehicles available $P_{u}$ is 7 ALS and 21 BLS. Belo Horizonte's EMS's response time standard is 10 minutes for ALS and 8 minutes for BLS. The maximum number of vehicles that each base can hold $C_{j}$ is two ambulances, according to the operational protocols stated. The response time between each potential ambulance location site $j$, and each demand node $i$ is taken using Google Maps.

The likely reliability of service for ambulance type $u \in U$ at demand point $i \in I$ at time period $t \in T$ given the allocation of $k \in K$ ambulances within its neighborhood $q_{iu k}^t$ is calculated as described on \eqref{rAD}. It is necessary to calculate the area-specific busy fraction of the $k$ servers within the same area surrounding $i$ ($b_{iuk}^t$). On its turn, $b_{iuk}^t$ is calculated as shown on \eqref{rAC}, corresponding to the time each demand site requires for each ambulance type $u$ each time period $t$, given the allocation of $k$ ambulances within a coverage zone.   The value of $k_{max}$ used is 3.

\subsection{Experimental scenarios}
The following scenarios are considered so that results can be analyzed and compared:

\begin{itemize}
    \item Scenario 1 (S1): based on base and ambulance configurations ruling in 2019. Thereby, variables $z_{j}$ and $x_{ju}^t$, which set base locations and ambulance allocations, are defined as parameters for both formulations, FLEET-ICt and LR-MEXCLP-ICt;
    \item Scenario 2 (S2): based on base configurations ruling in 2019 operations (variable $z_{j}$ is set as a parameter), both formulations work aiming to optimize the ambulance allocation considering one time period (static scenario);
    \item Scenario 3 (S3): based on base configurations ruling in 2019 operations (variable $z_{j}$ is set as a parameter), both formulations seek to optimize the ambulance allocation/relocation considering 24 time periods;
    \item Scenario 4 (S4): optimization for base location and ambulance allocation for both FLEET-ICt and LR-MEXCLP-ICt in a static approach (considering one time period only);
    \item Scenario 5 (S5): optimization for base location and ambulance allocation/relocation for both formulations considering 24 periods.
\end{itemize}

Belo Horizonte's EMS has a complex structure on the core of its operations that works as an operations center, a warehouse, and an ambulances base. For this reason, this operation center is set as an installed base for all five scenarios. Furthermore, for all experiments, the number of demands (26,024) and the number of ambulances (7 ALS and 21 BLS) are the same.

Both formulations are implemented using C++ and ran by IBM ILOG CPLEX solver on its 12.6.3 version. Experiments are realized on a computer HPE Server Proliant DL380 Gen9 with a processor Intel Xeon E5-2630 v4 (2.2GHz/10-core/25MB/85W) with 32GB Dual Rank x4 DDR4-240 memory.

\subsection{Results for the static and multi-period approaches of the deterministic mono-objective  FLEET-ICt }

In FLEET-ICt results, the optimal objective function value is the maximum number of calls covered according to the optimal EMS configuration. Using the number of covered calls for each EMS configuration, the coverage rate, which is expressed on percentage, can be computed as the ratio of the objective function value (number of demands covered) to the total demand. Further, in the same way,  the ALS and  BLS service coverage rate can be computed employing the number of demands covered for each type of ambulance. The results for the five scenarios proposed are presented in Table \ref{tab:20}.

\begin{center}
\begin{table}[h]
  \centering
  \caption{FLEET-ICt results}
    \begin{tabular}{lcccccc}
    \hline
    \multicolumn{1}{c}{\textbf{Results}} & \textbf{S1} & \textbf{S2} & \textbf{S3} & \textbf{S4} & \textbf{S5} \\
    \hline
    Stations installed: & 25    & 25    & 25    &  28    &  226\\
    Computational time: & --    & 0.07    & 19.33 & 29.48 & 2865.32\\
    Objective function value: & 15775 & 18220 & 18355 & 23383 & 23827\\
    Coverage rate ALS & 52.8\%  & 68.4\%  & 69.7\%  & 77.6\% & 80.3\%\\
    Coverage rate BLS & 65.9\%  & 71.1\%  & 71.1\%  & 98.2\%  & 99.2\%\\
    \hline
    \textbf{Total coverage rate} & \textbf{60.6\%} & \textbf{70.0\%} & \textbf{70.5\%} & \textbf{89.9\%} & \textbf{91.6\%}\\
    \hline
    \end{tabular}%
  \label{tab:20}%
\end{table}%
\end{center}

It can be seen that even considering the base locations ruling in 2019, the results from FLEET-ICt formulation bring advantages. In both approaches, static (S2) and its multi-period (S3), the results increase ALS and BLS coverage and the total coverage rate. These results point out that it is possible to improve the system performance by relocating ambulances among the periods.

Likewise, the comparison between scenarios S4 and S5, where base location and ambulance allocation/relocation are both optimized, shows that the objective function value and the coverage rate can be increased through the relocation of ambulances between the multi-period of the day (S5). The scenario S5 is the one with the highest objective function value (23,827) and the highest coverage rate. It is important to note that despite the difference between the static version's coverage rate, considering a single period, and the multi-period version being equal to 0.7\%, this difference means an increase of 444 calls covered by the multi-period approach.

There is no run time for scenario S1 because no optimization is performed. The location of stations and ambulance allocations are the previous set. The optimization run at scenario S2 spends 0.07 seconds. This scenario is the static case for one time period only, with the location of stations predetermined. The scenario S3 also considers predefined stations location and optimize the allocation and relocation of ambulances among stations in different periods, which gives the problem a dynamic aspect.  The computer takes 19.33 seconds to provide the results for the scenario S3 configuration. It takes 29.48 seconds to solve the problem according to scenario S4. This scenario considers station location and ambulance allocation statically, one time period only. Finally, the scenario S5 takes 2865.32 seconds to optimize the station location and ambulance allocation/relocation in a multi-period context.

According to ruling configurations in 2019 considered in this paper, the scenarios S1, S2, and S3 optimization establish that 24 stations have to be installed. In scenarios S4 and S5, the results propose 28 and 226 stations installed, respectively. The high number of stations in scenario S5 might be understood as a consequence of maximizing coverage in a multi-period approach. Since 28 ambulances are available for 24 periods, the optimization might suggest installing up to 672 bases ($28 \times 24$).


\subsection{Results for the static and multi-period approaches of the probabilistic mono-objective LR-MEXCLP-ICt}

In LR-MEXCLP-ICt, the objective function value is the number of covered demands according to the service level provided. Thus, the coverage rates are computed through the objective function. It is the ratio of the number of the covered demand, according to the service level provided, to the total demands. The results for the five proposed scenarios implemented through the LR-MEXCLP-ICt formulation are presented in Table \ref{tab:40}.

Considering the scenarios based on configurations ruling in 2019 (S1, S2, and S3), the lack of flexibility in the location of the bases results in a better objective function value and, consequently, the best total rate of coverage under the conditions established by scenario S2, when the problem is treated as static, without reallocation.  Furthermore, comparing the coverage rate for ALS and BLS,  the scenario S3 gets a better coverage rate for ALS, but scenario S2 attains a better coverage rate for BLS.

The results from S5 is better than S4. It suggests that the objective function value and the total coverage rate can be increased through ambulance relocation. However, it is observed that in scenario S4, the static allocation settings might generate a better coverage rate for BLS than that one provided by scenario S5.  Furthermore, the S5 scenario takes significantly more computational time than the other scenarios to find the optimal solution since it deals with location and allocation decisions and relocation decisions between periods. Consequently, in this scenario, more stations are installed.

\begin{center}
\begin{table}[h]
  \centering
  \caption{LR-MEXCLP-ICt results}
    \begin{tabular}{lcccccc}
    \hline
    \multicolumn{1}{c}{\textbf{Results}} & \textbf{S1} & \textbf{S2} & \textbf{S3} & \textbf{S4} & \textbf{S5}\\
    \hline
    Stations installed: & 25    & 25    & 25    & 27    &  179\\
    Computational time: & -     & 0.11     & 404.27 & 4.87 & 1637.02\\
    Objective function value: & 12855.6 & 15011.4 & 14742.1 & 18607.6 & 18757.2\\
    Coverage rate ALS & 52.7\%  & 68.3\%  & 69.4\%  & 77.2\% & 79.4\%\\
    Coverage rate BLS & 63.2\%  & 68.7\%  & 67.4\%  & 87.8\%  & 86.5\%\\
    \hline
    \textbf{Total coverage rate} & \textbf{59.0\%} & \textbf{68.5\%} & \textbf{68.2\%} & \textbf{83.5\%} & \textbf{83.6\%}\\
    \hline
    \end{tabular}%
  \label{tab:40}%
\end{table}%
\end{center}

Comparing the LR-MEXCLP-ICt and FLEET-ICt versions, note that the FLEET approach's total coverage rate overestimates the coverage rate obtained by the LR-MEXCLP-ICt approach. This result is expected because while the first approach considers allocating a single ambulance in a demand point's coverage area, it covers all calls generated in this area. The second approach considers that when the ambulance occupancy rates are very high, in periods of high demand, the allocation of a single ambulance in the coverage area of a demand point may not cover all calls.

\subsection{Results for the bi-objective approach}

The set of inequalities \eqref{rA14} on FLEET-ICt formulation and the set of inequalities \eqref{rA39} on LR-MEXCLP-ICt formulation limit that Belo Horizonte's EMS must allocate no more than two ambulances per station. Consequently, at least 14 stations have to be installed to meet these constraints. From this, the initial value defined for $\epsilon$ is 14. Furthermore, the maximum $\epsilon$ value is the maximum number of stations previously obtained, i.e. 226 for the FLEET-ICt and 179 for LR-MEXCLP-ICt. Hence $\epsilon$ values used in the $\epsilon$-constraint method are all integer from the minimum to the maximum value.

\begin{figure}[h]
\caption{Pareto set for the FLEET-ICt formulation - objective functions $F_1$ and $F_3$}
		\includegraphics[scale=0.55]{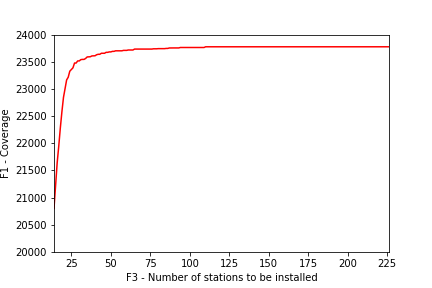}
\label{fig:MO_FL}%
\centering
\end{figure}

\begin{figure}[h]
\caption{Pareto set for the LR-MEXCLP-ICt formulation - objective functions $F_2$ and $F_3$}
		\includegraphics[scale=0.55]{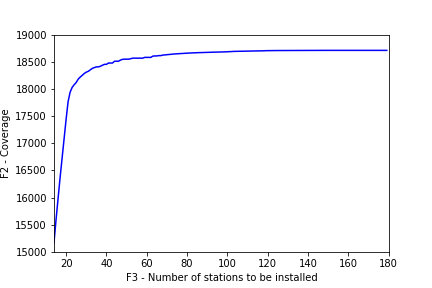}
\label{fig:MO_LR}%
\centering
\end{figure}

Figures \ref{fig:MO_FL} and \ref{fig:MO_LR} presents the value for both objective functions for different values of $\epsilon$ for the FLEET-ICt and LR-MEXCLP-ICt models, respectively. This figures show that to increase the number of covered demands, the number of stations installed also has to be increased. Since the second objective function aims to minimize the number of stations installed, it is impossible to improve the number of covered demands without compromising it.

When the number of stations is greater or equal to 28 in the multi-objective FLEET-ICt approach, the objective function value maximizes the number of covered demands starts to get relatively smaller increments. For 28 stations, the objective function value is equal to 23,476, while for 226 stations, the optimal solution has an objective function value of 23,775, which means an increment of less than 1.3\%. For the multi-objective LR-MEXCLP-ICt model, when the number of stations is superior to 30, the objective function value that maximizes the number of covered demands starts to have relatively smaller increments. From the solution for $\epsilon$ equal to 30 to $\epsilon$ equal to 179 stations, there is an increment of 2.45\% in the coverage.

\section{Conclusions}

The decision-making process is fundamental in the EMS context in both rich and emerging nations. Given the limited available resources, policymakers need effective methods for planning, prioritization, and decision-making. The formulations present in this paper are solved in a multi-period and multi-objective approach and seek to answer questions such as $i)$ How many stations should be installed? $ii)$ Where should each station be located? $iii)$ Which station should cover each demand site each period? $iv)$ How should ambulances be allocated on stations each period?

Results indicate improvements in coverage levels through the use of both formulations. FLEET-ICt formulation results show an increment up to 31\% on coverage levels. Comparing static and multi-period optimization scenarios for FLEET-ICt model, coverage levels can be increased by up to 21\% when ambulances are allowed to switch bases. LR-MEXCLP-ICt suggests improvements up to 24\% on coverage levels. Coverage levels can be increased by up to 15\% when ambulances are allowed to change their station through time periods at static and multi-period optimization comparison of LR-MEXCLP-ICt. However, finding solutions with greater coverage levels requires higher computational efforts due to the larger size of the corresponding mathematical models, which have a more significant number of variables and constraints than the static formulations.

The FLEET-ICt model implies that ambulances are always available whenever a request arises. It can be seen as a drawback; once ambulances become unavailable on any occasion, they are assigned to occurrences. Therefore, this might lead to an overestimation of the system's capacity once demand points with many incidents can not be covered sometimes. For example, the model might choose to install just one vehicle to cover the highest number of demands as possible (within a coverage zone), despite the intensity of demands within this same coverage zone.

The LR-MEXCLP-ICt model adds a reliability factor inside the objective function to guarantee that more calls will be attended within the desired response time. This factor may lead to the allocation of a high number of ambulances covering more significant demand points. However, LR-MEXCLP-ICt also has its drawbacks. It is crucial to notice that the concept of busy fraction for each demand site assumes that ambulances are only allowed to answer demands within its coverage zone. For example, they can only be considered busy when assigned to occurrences within their coverage zone. It does not accurately represent real-life since ambulances can be assigned to answer calls at any demand site, even outside their coverage zone (as long as they are the best option). However, such a drawback is familiar to most of the coverage models, and it is also a drawback to FLEET-ICt formulation.

The multi-period results for both models required a high number of stations to be installed. In this case, the proposed multi-objective approach is an excellent option because it tries to minimize the stations' number while maximizing the coverage levels. The decision-maker can choose a solution from the optimum solutions set represent by the Pareto solutions. 
Multi-objective results suggest that only modest increments in the objective function are possible when the stations' number exceeds 28 for the FLEET-ICt formulation. In contrast, for the LR-MEXCLP-ICT model, increments in the objective function start to narrow below 30 installed stations.

Future studies might consider other objective functions, according to the needs that concern the system's decision-making process. For example, upper bounds on the number of relocation or the traveled distance might also be considered, evaluating the trade-off between service level and relocation costs. Another interesting future study is to compare the proposed models' solutions using a simulation model.

\section*{Acknowledgments}

This research is supported by SAMU-BH, PPSUS, and the Brazilian funding agencies CNPq, CAPES, and Fapemig. We, with this, sincerely thank these organizations for their support.

\biboptions{round}
\bibliography{mybibfile}

\end{document}